\documentclass[12pt]{article}
\usepackage[cp1251]{inputenc}
\usepackage[english]{babel}
\usepackage{amsfonts}
\usepackage{amssymb}
\usepackage{color}
\usepackage{epic}
\usepackage{graphics}
\usepackage{graphicx}

\newtheorem{Definition}{\indent Definition}[section]

\newtheorem{Lemma}{\indent Lemma}[section]

\begin{document}

\begin{center}
{\bf Fractional Calculus in Russia at the end of XIX century}
\end{center}

\vspace{5mm}
\begin{center}
{\bf Sergei Rogosin$^a$ and Maryna Dubatovskaya$^a$}
\end{center}

\begin{center}
$^a${\bf Department of Economics, Belarusian State University, 4,
Nezavisimosti ave, Minsk, 220030, Belarus; E-Mail:
rogosinsv@gmail.com}
\end{center}


\bigskip
\begin{abstract}
{In this survey paper we analyze the development of Fractional Calculus in Russia
at the end of XIX century, in particular, the results by A.V.Letnikov, N.Ya.Sonine and P.A.Nekrasov. Some of the discussed results are either unknown or inaccessible.}
\\[1mm]
{\noindent{\it Mathematics Subject Classification:} {primary  26A33, secondary 34A08, 34K37, 35R11, 39A70}\\[1mm]
}
\end{abstract}





\section{Introduction}

The year of the birth of Fractional Calculus is considered 1695 when Leibniz discussed the possibility of introducing the derivative of an arbitrary order in his letters to Wallis and Bernulli. Several attempts were made to give a precise meaning to this new notion. A comprehensive detailed analysis of the history of Fractional Calculus is given in \cite{SKM}. One of most productive periods in this history was middle-end of XIX century. Here we can mention works by Legendre, Fourier, Peacock, Kelland, Tardi, Roberts and others. The most advanced approach to the determination of the fractional derivative of an arbitrary order was proposed by Liouville. A deep analysis of the results on this subject was given in the article \cite{Let68a} by Letnikov. In particular, he recognized the basic role of the Liouville'as approach. Letnikov said \cite[p. 92] {Let68a}: "...we give a survey of the results by Liouville whom we ought to consider as the first scientist paid a serious attention to clarifying the question on the derivative of an arbitrary order. In 1832 he started to publish a series of articles devoted to the foundation and the application of his theory of general differentiation which is the first complete discussion of this topic. Before this work  were made only few but very important although not completely clear remarks on this subject."

It should be noted that the works by A.V.Letnikov constitutes the first rigorous and comprehensive construction of the theory of integro-differentiation. An extended description of the results by Letnikov is presented in the articles \cite{Pot03,Pot12} and in the book \cite{LetChe11} written in Russian.

In the middle-end of XIX century an interest to Fractional Calculus in Russia grew essentially \cite{Let68,Let73,Son72,Nek88,Vas-Zak61}. One of the reasons for it was a high standard in the research in Real and Complex Analysis in this period. Russian Universities took care of the level of the education of young scientists. For many professorship expectant there were given possibility to spend one-two years at the leading research centers and attend the lectures of known mathematicians.

Letnikov's results attracted people to this branch of the Science at least in Russia. Anyway these works remained unknown abroad and for a long time was unaccessible. After contribution by Letnikov the serious works on Fractional Calculus in Russia at the second part of XIX century were performed by N.Ya.Sonine and P.A.Nekrasov. They introduced  the complex-analytic technique into the study and application of derivatives and integrals of an arbitrary order.It should be noted that Complex Analysis was traditionally highly developed discipline starting from Leonard Euler who worked for a long period in Russia (1726-1741 and 1776-1783). This part of Mathematical Analysis was essentially developed in XIX century by M.V.Ostrogradsky, V.Ya.Bunyakovsky, P.L.Chebyshev, A.M.Lyapunov and many others. In particular, Sonine and Nekrasov found a fractional analog of the classical Cauchy integral formula for analytic functions.

In our article we describe the contribution of Alexey Vasil'evich Letnikov (1837-1888), Nikolai Yakovlevich Sonine (1849-1915) and Pavel Alekseevich Nekrasov (1853-1924) to Fractional Calculus and the role of these results in the modern Fractional Calculus and its Applications.

\section{ Liouville approach and its analysis by Letnikov}

As it was already said, A.V.Letnikov considered (see, e.g. \cite{Let68,Let68a,Let73}) that the Liouville's theory constitutes the only complete treatment of differentiation of an arbitrary order. Realizing the great importance of this theory, Lentnikov had seen that its certain parts did not receive a proper justification and led to some misunderstanding in the works of Liouville's followers.

Let us present here Letnikov's description of the elements of the Lioville's theory following \cite{Let68a}. He started his analysis with the  definitions given by Liouville.
\begin{Definition}
\label{Lio_def}
Let the function $y(x)$ is represented in the form of the following series of exponents
\begin{equation}
\label{Lio_def1}
y(x) = A_1 e^{m_1 x} + A_2 e^{m_2 x} + \ldots,
\end{equation}
which is denoted for shortness as $\sum A_m e^{m x}$.

Fractional derivative of the order $p$ is defined by multiplying each term of the series by $p$-th power of the index $m$:
\begin{equation}
\label{Lio_def2}
\frac{d^p y}{d x^p} = \sum A_m m^p e^{m x}.
\end{equation}

If $p$ is negative then formula (\ref{Lio_def2}) determined the fractional integral of order $-p$.
\end{Definition}
Fractional integral of order $-p$ is denoted by Liouville as $\int{}^{-p} y d x^{-p}$. Liouville considered this definition as the only possible way to generalize the usual derivative. Evaluating its role, Letnikov stressed that Definition \ref{Lio_def} contains a key ideas to establish a deep analogy with differences and powers and thus could lead to a more simple construction.

Anyway, the above definition had a very important restriction. It cannot be applied to an arbitrary function since not all of them possess representations in series of exponents. Liouville himself understood this difficulty. He proposed a way to overcome it. By performing the change of variable $z = e^x$ for the function $y = F(x)$ one can expand the composite function $y = F(\ln z)$ (with $x = \ln z$) in a converging power series
\begin{equation}
\label{ser1}
F(\ln z) = \sum A_m z^m.
\end{equation}
Thus, the initial function $y = F(x)$ admits representation via series of exponents
\begin{equation}
\label{ser2}
F(x) = \sum A_m e^{m x}.
\end{equation}
But possibility to represent $y = F(\ln z)$ in form (\ref{ser1}) met several restrictions. For instance, if we suppose to get representation of $y = F(\ln z)$ in a form of series in positive powers of $z$, then all derivatives of   $y = F(x)$ at $x = \infty$ should be equal to zero since
$$
\left(F(\ln z)\right)^{\prime}_{z} = \frac{F^{\prime}(\ln z)}{z}, \left(F(\ln z)\right)^{\prime\prime}_{z} = \frac{F^{\prime\prime}(\ln z) - F^{\prime}(\ln z)}{z^2}, \ldots
$$
Similar restriction appears if we suppose to represent  $y = F(\ln z)$  in a form of series in negative powers of $z$ since we deal in this case with the function $y = F(- \ln z)$. Such conditions looks fairly strong. Moreover they are neither necessary nor sufficient for the representation of the type (\ref{ser2}).

Liouville met such a restriction trying to calculate the derivative of a fractional order of the  power function. He started from the Euler formula
$$
\frac{1}{x^m} = \frac{\int\limits_{0}^{\infty} e^{- x z} z^{m - 1} dz}{\Gamma(m)}.
$$
Liouville supposed that the above integral can be represented in a form of the exponential sum $\sum A_n e^{-nx}$. Here all coefficients in such representation should be infinitely small. Then using his main definition Liouville arrived at the formula of the derivative of this function
\begin{equation}
\label{FD_powers}
\frac{d^p}{d x^p} \frac{1}{x^m} = \frac{\int\limits_{0}^{\infty} e^{- x z} (- z)^{p} z^{m - 1} dz}{\Gamma(m)}.
\end{equation}
Thus by definition of $\Gamma$-function we get after substitution $xz = t$ the following formula
\begin{equation}
\label{FD_powers1}
\frac{d^p}{d x^p} \frac{1}{x^m} = \frac{(-1)^{p} \Gamma(m+p)}{\Gamma(m) x^{m+p}}.
\end{equation}
In his first articles Liouville used the only definition of the $\Gamma$-function of positive variable (later he noted that he was not familiar with the general definition of the $\Gamma$-function by Legendre and Gauss). Therefore he supposed that due to assumptions $m > 0$, $m + p > 0$ one needs to use in the above definition so called auxiliary functions. Being very important, the use of auxiliary functions did not lead to a general definition of the fractional derivative. Liouville showed that if one supposed an existence of auxiliary functions then these necessarily had to be entire functions.

Letnikov claimed and proved that it follows from his analysis that Liouville's formulas were so general that had no need of any auxiliary function. Later on, several attempts to correct Liouville's approach were made. In particular, Letnikov analyzed in \cite{Let68a} the works by Kelland, Tardi and Roberts. But really rigorous approach which transformed Liouville's constructions to the general definition of the fractional derivative was proposed by Letnikov.

\section{Letnikov's contribution to Fractional Calculus}
\label{Letnikov}

\subsection{ Letnikov or Gr\"unwald-Letnikov derivative}
\label{LetnikovFD}

Starting his work on determination of the derivative of an arbitrary order, Letnikov posed this problem (\cite{Let68,Let73}) as interpolation in form of the elements of two sequences consisting of successive derivatives of the function $f(x)$
\begin{equation}
\label{seq_der}
(a)\;\;\;\;\; f(x), f^{\prime}(x), f^{\prime\prime}(x), \ldots, f^{(n)}(x), \ldots
\end{equation}
and of successive $n$-fold integrals of this function
\begin{equation}
\label{seq_int}
(b)\;\;\;\;\; f(x), \int f(x) dx, \int^{2} f(x) d x^2, \ldots, \int^{n} f(x) d x^n, \ldots
\end{equation}
In other words he tries to find such a formula of the derivative of an arbitrary order $\alpha$ which for nonegative integer
$\alpha = 0, 1, 2, \ldots$ coincides with the corresponding elements of the sequence $(a)$ and for nonpositive integer $\alpha = 0, - 1, - 2, \ldots$ coincides with the corresponding elements of the sequence $(b)$. Denoting this formula by
$$
D^{\alpha} f(x)\;\;\; {\mathrm{or}} \;\;\; \frac{d^{\alpha} f(x)}{d x^{\alpha}}
$$
he try to get this new object having (whenever it is possible) same properties as elements of sequence $(a)$ or $(b)$ when $\alpha$ is integer.

Next idea by Letnikov was to restrict the generality of the above question and to consider instead of the sequence $(b)$ (of indefinite $n$-fold integrals) the sequence of definite integrals supposing that $f(x)$ is continuous on certain interval $[a, x]$, i.e. to interpolate in form elements of the double sequence
\begin{equation}
\label{seq_int-der}
(A)\;\;\;\;\; \ldots \int\limits_{a}^{x}\int\limits_{a}^{x}f(x) dx^2, \int\limits_{a}^{x}f(x) dx, f^{\prime}(x), f^{\prime\prime}(x), \ldots,
\end{equation}
in which any element is the derivative of the previous one.

The corresponding interpolating object he denoted
$$
\left[D^{\alpha} f(x)\right]_{a}^{x}.
$$
In order to get such interpolation Letnikov proposed to examine the following formula
\begin{equation}
\label{Let_def1} 
\frac{\sum\limits_{k=0}^{n} (-1)^k \left(\begin{array}{c} \alpha \\
k \end{array}\right) y(x - k h)}{h^{\alpha}},
\end{equation}
where $h = \frac{x - a}{n}$ and $\left(\begin{array}{c} \alpha \\
k \end{array}\right)$ denotes the binomial coefficient. This approach was independently used by Gr\"unwald \cite{Gru67} and by Letnikov \cite{Let68}. When Letnikov found the paper by Gr\"unwald, he decided to decline publication of his work, but later changed his mind. Letnikov developed in \cite{Let68} more rigorously than in \cite{Gru67} the theory of the derivative of an arbitrary order and found its relationship with many results known in this area.

Elementary algebra yields that, for $\alpha = m$ being positive integer number, the derivatives of the corresponding order can be defined as a limit of the above expression
\begin{equation}
\label{der_alpha=m}
f^{(m)}(x) = \lim\limits_{\delta\rightarrow 0}
\frac{f(x) - \left(\begin{array}{c} m \\ 1\end{array}\right) f(x - \delta) + \left(\begin{array}{c} m \\ 2\end{array}\right) f(x - 2 \delta) + \ldots + (-1)^n \left(\begin{array}{c} m \\ n\end{array}\right) f(x - n \delta)}{\delta^m}.
\end{equation}
Here ${\delta\rightarrow 0}$ is equivalent to ${n\rightarrow \infty}$, but the sum in the numerator remains finite since all binomial coefficients with $n > m$ vanishing. Thus, formula (\ref{der_alpha=m}) can be taken as the definition of the derivative of order $m\in {\mathbb N}$.

Vice versa, for $\alpha = - m$ being negative integer the expression under the limit sign in the right hand-side of right (\ref{der_alpha=m}) equals to
\begin{equation}
\label{int_alpha}
\frac{f(x) + \left(\begin{array}{c} m \\ 1\end{array}\right) f(x - \delta) + \left(\begin{array}{c} m \\ 2\end{array}\right) f(x - 2 \delta) + \ldots + \left(\begin{array}{c} m \\ n\end{array}\right) f(x - n \delta)}{\delta^{-m}}.
\end{equation}
Letnikov showed \cite[P. 5-12]{Let68} that the limit of this expression as ${\delta\rightarrow 0}$ or equivalently as ${n\rightarrow \infty}$ is equal to the multiply integral, i.e. (in his notation)
\begin{equation}
\label{int_alpha=-m}
\left[D^{-m} f(x)\right]_{a}^{x} = \lim\limits_{\delta\rightarrow 0}
\frac{f(x) + \left(\begin{array}{c} m \\ 1\end{array}\right) f(x - \delta) + \left(\begin{array}{c} m \\ 2\end{array}\right) f(x - 2 \delta) + \ldots + \left(\begin{array}{c} m \\ n\end{array}\right) f(x - n \delta)}{\delta^{-m}} =
\end{equation}
$$
= \int\limits_{a}^{x} d x_1 \int\limits_{a}^{x_1} d x_2 \ldots \int\limits_{a}^{x_{m-1}} f(x_m) d x_m.
$$
This magnitude $\left[D^{-m} f(x)\right]_{a}^{x}$ satisfies certain properties. First of all, if we apply to it similar operation of order $- p$, $p > 0$, then we will have
$$
\left[D^{-p} D^{-m} f(x)\right]_{a}^{x} = \left[D^{-m-p} f(x)\right]_{a}^{x}.
$$
Next, if we take the derivative $\frac{d^p}{d x^p}$ of order $p > 0$, then we will have
$$
\frac{d^p}{d x^p} \left[D^{-m} f(x)\right]_{a}^{x} = \left[D^{-m+p} f(x)\right]_{a}^{x}, \;\;\; if \;\;\; m > p,
$$
and
$$
\frac{d^p}{d x^p} \left[D^{-m} f(x)\right]_{a}^{x} = \frac{d^{p-m} f(x)}{d x^{p-m}}, \;\;\; if \;\;\; m < p.
$$
Thus, in particular, the symbol $\left[D^{-m} f(x)\right]_{a}^{x}$ means $m$-times differentiable function whose all derivatives up to $m$-th order are vanishing at $x= a$.

Formulas (\ref{int_alpha=-m}) and (\ref{der_alpha=m}) coincide with the corresponding elements of the double sequence $(A)$. Therefore, it leads Letnikov to the conclusion that the following limit
\begin{equation}
\label{GruLet0}
[D^{\alpha}f(x)]_{a}^{x} := \lim\limits_{\delta\rightarrow 0}
\frac{f(x) - \left(\begin{array}{c} \alpha \\ 1\end{array}\right) f(x - \delta) + \left(\begin{array}{c} \alpha \\ 2\end{array}\right) f(x - 2 \delta) + \ldots + (-1)^n \left(\begin{array}{c} \alpha \\ n\end{array}\right) f(x - n \delta)}{\delta^{\alpha}}
\end{equation}
is a good candidate to solve the interpolation problem for the sequence $(A)$, i.e. to be the  derivative of arbitrary order.

The relations of this new object to known formulas of the fractional derivatives (integrals) were described by Letnikov  \cite[P. 15]{Let68} using the following elementary
\begin{Lemma}
\label{summation}
Let $(\alpha_k)$ be a sequence of (real or complex) numbers such that
$$
\lim\limits_{k\rightarrow \infty} \alpha_k = 0\;\;\; and \;\;\; \lim\limits_{k\rightarrow \infty} (\alpha_1 + \alpha_2 + \ldots + \alpha_k) = C,
$$
and $(\beta_k)$ be a sequence of (real or complex) numbers such that
$$
\lim\limits_{k\rightarrow \infty} \beta_k = 1.
$$

Then the sequence of their products has the limits equal to $C$, i.e.
$$
\lim\limits_{k\rightarrow \infty} (\alpha_1 \beta_1 + \alpha_2 \beta_2 + \ldots + \alpha_k \beta_k) =
\lim\limits_{k\rightarrow \infty} (\alpha_1 + \alpha_2 + \ldots + \alpha_k) = C.
$$
\end{Lemma}

The above formula is valid for $[D^{\alpha}f(x)]_{a}^{x}$ with $\alpha < 0$ (i.e. for fractional integral of order $-\alpha$ in modern language).

The corresponding justification of the formula of $[D^{\alpha}f(x)]_{a}^{x}$ with $\alpha > 0$
(i.e. representation of fractional derivative) Letnikov supposed additional that the function $f(x)$ is $(n+1)$-times continuously differentiable on the interval $(a, x)$, where $n$ is a largest integer smaller than $\alpha$, i.e. $n < \alpha < n + 1$. Then using quite cumbersome transformation of the binomial coefficients \cite[P. 21-26]{Let68} he has got that the limit in (\ref{GruLet0}) is equal
\begin{equation}
\label{GruLet0a}
[D^{\alpha}f(x)]_{a}^{x} = \frac{f(a) (x - a)^{-\alpha}}{\Gamma(-\alpha + 1)} + \frac{f^{\prime}(a) (x - a)^{-\alpha + 1}}{\Gamma(-\alpha + 2)} + \ldots + \frac{f^{(n)}(a) (x - a)^{-\alpha + n}}{\Gamma(-\alpha + n + 1)} +
\end{equation}
$$
 + \frac{1}{\Gamma(-\alpha + n + 1)} \int\limits_{a}^{x} (x - \tau)^{- \alpha + n} f^{(n + 1)}(\tau) d \tau.
$$
Note that the same result is true if $\alpha\in {\mathbb C}$, ${\mathrm{Re}} \alpha > 0$. Integration by parts showed that (\ref{GruLet0a}) can be taken as the definition of fractional derivative of arbitrary order $\alpha > 0$. Slightly more general form can be written for any $s\in {\mathbb Z}$, $s \geq n$, $n < {\mathrm{Re}} \alpha < n + 1$ (of course under additional smoothness conditions)
\begin{equation}
\label{GruLet_gen}
[D^{\alpha}f(x)]_{a}^{x} = \sum\limits_{k=0}^{s} \frac{f^{(k)}(a) (x - a)^{-\alpha + k}}{\Gamma(-\alpha + k + 1)} + \frac{1}{\Gamma(-\alpha + s + 1)} \int\limits_{a}^{x} (x - \tau)^{- \alpha + s} f^{(s + 1)}(\tau) d \tau.
\end{equation}

In \cite{Let68} Letnikov paid attention to relationship of his formulas with known constructions. In particular, he showed that if the function $f(x)$ is defined, infinitely differentiable on $[x, \infty)$ and vanishes together with any derivative when $x$ is tending to $\infty$, then the following formula hold for any $\alpha, {\mathrm{Re}}\, \alpha < 0$,
$$
[D^{\alpha}f(x)]_{+\infty}^{x} = \frac{1}{\Gamma(-\alpha)} \int\limits_{+\infty}^{x} (x - \tau)^{-\alpha - 1} f(\tau) d \tau = \frac{1}{(-1)^{\alpha} \Gamma(-\alpha)} \int\limits_{0}^{+\infty} z^{-\alpha - 1} f(x + z) d z,
$$
i.e. coincides with the corresponding integral defined by Liouville. Similarly, for any $\alpha, 0 \leq n < {\mathrm{Re}}\, \alpha < n + 1$ Letnikov discovered that
$$
[D^{\alpha}f(x)]_{+\infty}^{x} = \frac{1}{\Gamma(-\alpha + n + 1)} \int\limits_{+\infty}^{x} \frac{f^{(n + 1)}(\tau) d \tau}{(x - \tau)^{\alpha-n}} =
\frac{1}{(-1)^{\alpha - n - 1} \Gamma(-\alpha + n + 1)} \int\limits_{0}^{+\infty}  \frac{f^{(n + 1)}(x + z) d z}{z^{\alpha - n}}.
$$
He also noted that the considered class of functions is not empty, it contains, in particular, all functions of the form $x^{m} e^{- x}$.

In \cite{Let68} Letnikov also presented a series of formulas for the values of his derivative of an arbitrary order of elementary functions such as power function $(x - a)^{\beta}$, exponential function $e^{m x}$, logarithmic function $\log x$, exponential-trigonometric functions $e^{\beta x} \sin \gamma x$, $e^{\beta x} \cos \gamma x$, rational functions $\frac{P(x)}{Q(x)}$. These formulas coincide with nowadays known formulas (see, e.g. \cite{SKM,GKMR2}). Composition formulas for fractional derivatives and integrals were found in \cite{Let68} too. Last result, presented in \cite{Let68} was so called Leibniz rule for the fractional derivative/integral of the product of functions. Note, that after the death of A.V. Letnikov it was created a committee examined some of his manuscripts \cite{NekPok89}. A few results were then published, but not all were found. In particular, the members of the committee reported that they didnot find any results on Abel integrals as it was expected some researchers.

\subsection{Solution to certain differential equations}
\label{LetnikovDE}

In \cite{Lio32a} Liouville made a background for further development of Fractional Calculus. In order to show an importance of
the new branch of Science he solved in \cite{Lio32b} a number of problems (mainly from geometry, classical mechanics and mathematical physics)
by using his constructions of integral and derivatives of an arbitrary order. Later in \cite{Lio34} he also discussed the tautochrone problem and usage of fractional derivatives to its solution.

In his master thesis  Letnikov carefully examine these results by Liouville and came to the conclusion that Liouville's solutions of the problems can be obtained by using more traditional methods too. He also remarked that incorrect using of Liouville construction by his followers led to certain misunderstanding and even mistakes. Note that master thesis by Letnikov was reprinted in Russian recently in \cite{LetChe11,Pot12}.

 Anyway Letnikov believed that newly created technique could find proper applications. One of these applications was presented in his article \cite{Let89a} devoted to the use of fractional derivative to the solution of the differential equation
\begin{equation}
\label{DE}
(a_n + b_n x) \frac{d^n y}{d x^n} + (a_{n-1} + b_{n-1} x) \frac{d^{n-1} y}{d x^{n-1}} + \ldots + (a_0 + b_0 x) y = 0.
\end{equation}
These results were lectured by Letnikov at the meeting of Mathematical Society on April 16th, 1876, and at the meeting of the Warsaw Congress of naturalists on September 3rd, 1876. They were reprinted by P.A.Nekrasov who parsed the Letnikov's archive after his death.

Denoting
$$
\varphi(\rho) := a_n \rho^n + a_{n-1} \rho^{n-1} + \ldots + a_1 \rho + a_0, \; \psi(\rho) := b_n \rho^n + b_{n-1} \rho^{n-1} + \ldots + b_1 \rho + b_0,
$$
equation (\ref{DE}) can be rewritten in the following symbolic form
\begin{equation}
\label{DE1}
\varphi\left(\frac{d}{d x}\right) y + x \psi\left(\frac{d}{d x}\right) y = 0.
\end{equation}

Suppose that equation
\begin{equation}
\label{ChE}
\psi(\lambda) = 0
\end{equation}
has different zeroes $\lambda_1, \lambda_2, \ldots, \lambda_n$. Denoting for each $j = 1, 2, \ldots, n$
\begin{equation}
\label{ChE1}
y := e^{\lambda_j x} Y
\end{equation}
one can rewrite (\ref{DE1})
\begin{equation}
\label{DE2}
\varphi\left(\lambda_j + \frac{d}{d x}\right) Y + x \psi\left(\lambda_j + \frac{d}{d x}\right) Y = 0.
\end{equation}
Crucial idea by Letnikov was to look for the solution to equation  (\ref{DE2})
in the form
\begin{equation}
\label{DE3}
Y = \left[D^p y_j\right]_{a}^{x},
\end{equation}
where $y_j$ is a new unknown function and $\left[D^p \right]_{a}^{x}$ is an interlimit (Letnikov-type) derivative whose order $p$ to be defined later.

Let there exists the function $y_j$ vanishing at $x = a$ together with all derivatives up to order $n - 2$ satisfying the following equation
\begin{equation}
\label{DE4}
\varphi_j\left(\frac{d}{d x}\right) y_j + x \psi_j\left(\frac{d}{d x}\right) y_j = 0,
\end{equation}
where
$$
\psi_j(\rho) = \frac{\psi(\lambda_j + \rho)}{\rho},\;\;\; p + 1 = \frac{\varphi(\lambda_j)}{\psi_j(0)} = \frac{\varphi(\lambda_j)}{\psi^{\prime}(\lambda_j)} = A_j,\;\;\; \varphi_j(\rho) = \frac{\varphi(\lambda_j + \rho) - A_j \psi_j(\rho)}{\rho}.
$$
Then there exists the solution to equation  (\ref{DE1}) which is represented in the form
\begin{equation}
\label{Sol_DE1}
y = e^{\lambda_j x} \left[D^{A_j - 1} y_j\right]_{a}^{x}.
\end{equation}
Thus, by this transformation we reduce equation  (\ref{DE}) of order $n$ to equation  (\ref{DE4}) of order $n - 1$. By applying this method one can reduce the order of the equation to 1 and get the possible solution via successive application of the interlimit derivative.

\section{Sonine-Letnikov discussion}

In 1868 A.V. Letnikov published the main part of his master thesis as an article in Mathematical Sbornik \cite{Let68} supplemented by the historical survey on the development of the theory of differentian of an arbitrary order \cite{Let68a}. These articles and Gr\"unwald's article \cite{Gru67} was criticized by N.Ya.Sonine in \cite{Son72} who also present in  \cite{Son72} his own approach to creation of the derivatives of an arbitrary order.
Sonine's article started with the discussion of Liouville's definition of the derivative of an arbitrary order (not necessarily positive). The latter definition is based on the derivative of an arbitrary order $p\in {\mathbb R}$ of exponential function
$$
\frac{d^p}{d x^p} e^{mx} = m^p e^{mx}
$$
and on the possibility to expand a differentiable function into series in exponents (the Dirichlet series in modern language):
$$
f(x) = \sum\limits_{k=-\infty}^{+\infty} A_k e^{m_k x}.
$$
Sonine made two important remarks concerning Liouville's definition. First, he showed that the derivative of negative order (i.e. the fractional integral) cannot be considered as inverse to the fractional derivative. Second remark by Sonine is related to the problem discovered by Liouvile himself. If one uses Liouville's definition of the derivative of arbitrary order from power function, then it leads immediately to a kind of contradiction. Liouville founded this phenomena using the following representation of $x$
$$
x = \lim\limits_{\beta\rightarrow 0} \frac{e^{\beta x} - e^{-\beta x}}{2 \beta}.
$$
If we suppose that the limit and the fractional derivative are interchangeable then the half-derivative of $x$ becomes infinite. Liouville concluded from it that there exist additional functions equal to the derivative of $f(x) = 0$. Contrariwise, Sonine has shown that such contradiction follows from not completely rigorous way of expansion of the function into the series in exponents. He also remarked that Liouville's proof of existence of additional functions does not have a proper rigor.

In the second part of his article \cite{Son72} Sonine criticized an approach by Gr\"unwald and Letnikov in which the fractional derivative is defined by the following limit
\begin{equation}
\label{GruLet1}
D^{\alpha}[f(x)]_{x=a}^{x=x} := \lim\limits_{\delta\rightarrow 0}
\frac{f(x) - \left(\begin{array}{c} \alpha \\ 1\end{array}\right) f(x - h) + \left(\begin{array}{c} \alpha \\ 2\end{array}\right) f(x - 2 h) + \ldots + (-1)^n \left(\begin{array}{c} \alpha \\ n\end{array}\right) f(x - n h)}{h^{\alpha}},
\end{equation}
where $n h = x - a$.

Sonine had two main objections. First, he noted that the series in the numerator of (\ref{GruLet1}) is converging. Hence the fraction should be infinite. Second remark by Sonine that if we apply to the fractional integral $D^{-\alpha}$ the fractional derivative $D^{\beta}$ ($\alpha,\beta > 0$), then by Leibnitz rule the result should coincide with $D^{\beta-\alpha}$. It leads to contradiction even for the function $f(x) = 1$ since it should exist for any $\beta$, but it is so only for ${\mathrm{Re}} \alpha > [{\mathrm{Re}} \beta]$.

Concerning the first remark, Letnikov noted in \cite{Let73} that in the case ${\mathrm{Re}} \alpha > 0$ formula (\ref{GruLet1}) determines the fractional integral (coinciding with m-times repeated integral if $\alpha = m\in {\mathbb N}$) if the series in the numerator of (\ref{GruLet1}) is converging and its sum is equal to zero. Moreover, Letnikov gave sufficient conditions for existence of the limit in (\ref{GruLet1}).

Second question by Sonine appeared due to his incorrect application of the Leibnitz rule. Letnikov noted that fractional integral and fractional derivative are defined by different formulas:
\begin{equation}
\label{GruLet_int}
D^{\alpha}[f(x)]_{x=a}^{x=x} = \frac{1}{\Gamma(-\alpha)} \int\limits_{a}^{x} (x - t)^{-\alpha - 1} f(t) dt,\;\;\; \alpha < 0,
\end{equation}
\begin{equation}
\label{GruLet_der}
D^{\alpha}[f(x)]_{x=a}^{x=x} = \sum\limits_{k=0}^{m} \frac{f^{(k)}(a) (x - a)^{-\alpha + k}}{\Gamma(-\alpha+ k + 1)} + \frac{1}{\Gamma(-\alpha + m + 1)} \int\limits_{a}^{x} (x - t)^{-\alpha + m} f^{(m+1)}(t) dt,
\end{equation}
$$
\alpha > 0, m = [\alpha].
$$
Both formulas are applied under certain conditions. Thus successive application of these two formulas can lead to certain contradiction if we do not take into account the above conditions.

\section{Sonine's contribution to Fractional Calculus}

\subsection{Sonine's fractional derivative and integral}

In his polemical article \cite{Son72} N.Ya.Sonine not only criticized  Gr\'unwald-Letnikov approach, but also proposed
another form of ``general'' fractional derivative. For his formula Sonine used generalization of the Cauchy integral (or better to say, Cauchy-type integral since below integral is defined for any continuous function $f$):
\begin{equation}
\label{Sonine1}
\frac{d^{\alpha} f(x)}{d x^\alpha} = \frac{\Gamma(p + 1)}{2\pi i} \int\limits_{\gamma} \frac{f(\tau) d \tau}{(\tau - x)^{\alpha + 1}},
\end{equation}
where $\gamma$ is a closed simple smooth curve on the complex plane surrounding the point $x$ (without loss of generality one can assume that  $\gamma$ is the circle of radius $r$ centered at the point $x$). This formula is really a good candidate for generalization of usual derivatives since for $\alpha = p$ positive integer (\ref{Sonine1}) gives the value of $p$-th derivative at the point $x$ assuming that the function $f$ is $p$-time differentiable at $x$.

This formula was analyzed by Letnikov in his answer \cite{Let73} on the remarks by Sonine. He proved that under assumption that the function $f$ is $(m + 1) = ([\alpha] + 1)$-times continuously differentiable inside the circle  $\gamma$ formula (\ref{Sonine1}) coincides with Letnikov's formula (\ref{GruLet_der}). {Note that Letnikov discussed Sonine's formula  (\ref{Sonine1}) under the stronger assumptions on the function $f$.}
\begin{equation}
\label{SonLet}
\frac{d^{\alpha} f(x)}{d x^\alpha} = \frac{(-r)^{-\alpha} f(x + r)}{\Gamma(-\alpha + 1)} + \frac{(-r)^{-\alpha+1} f^{\prime}(x + r)}{\Gamma(-\alpha + 2)} + \ldots + \frac{(-r)^{-\alpha+m} f^{(m)}(x + r)}{\Gamma(-\alpha + m + 1)} +
\end{equation}
$$
+ \frac{\Gamma(\alpha - m)}{2\pi r^{\alpha - m - 1}} \int\limits_{0}^{2\pi} f^{(m+1)}(x + r e^{i \theta}) e^{-(\alpha - m - 1) i \theta} d \theta.
$$
Since the last integral can be transformed to the form
$$
\frac{\Gamma(p - m)}{2\pi r^{\alpha - m - 1}} \int\limits_{0}^{2\pi} f^{(m+1)}(x + r e^{i \theta}) e^{-(\alpha - m - 1) i \theta} d \theta = \frac{(-1)^{-\alpha + m}}{\Gamma(-\alpha + m + 1)} \int\limits_{x + r}^{x} f^{(m+1)}(\tau) (\tau - x)^{-\alpha + m} d \tau,
$$
then formula (\ref{Sonine1}) coincides with the definition of the fractional derivative given by Letnikov  (\ref{GruLet_der}). {In \cite{Son72} Sonine concluded that his formula cannot coincide with Gr\"unwald-Letnikov formula for $\alpha > 0$ without adding auxiliary function.}

Sonine's definition of fractional derivative of negative order (i.e. fractional integral) has been criticized  by Letnikov in \cite{Let73}. Sonine used Leibniz rule (composition formula) for fractional derivatives (which does not valid in general, see \cite{SKM}). Fractional integral by Sonine is defined as inverse operation for fractional derivative
\begin{equation}
\label{Leibniz}
\frac{d^{\alpha}}{d x^\alpha} \frac{d^{-\alpha} f(x)}{d x^{-\alpha}} = f(x).
\end{equation}
From this formula Sonine concluded that it should be existed so called auxiliary function $\psi(x)$ satisfying the following relation
\begin{equation}
\label{aux1}
\frac{d^{\alpha} \psi(x)}{d x^\alpha} = 0.
\end{equation}
Sonine took the function $\psi(x)$ in the form
$$
\psi(x) = A_{1} (x - a)^{\alpha - 1} + A_{2} (x - a)^{\alpha - 2} + \ldots + A_{k} (x - a)^{\alpha - k},
$$
where $A_j, j=1, 2, \ldots, k,$ are arbitrary constants, $k = [p] + 1$, but $a$ is not defined by Sonine. If Cauchy formula is taken as the definition of the derivative of an arbitrary order, then the auxiliary function has to satisfy the relation
$$
\int\limits_{\gamma} \frac{\psi(\tau) d \tau}{(\tau - x)^{\alpha + 1}} = 0.
$$
Since it was shown that by integration by parts  formula (\ref{Sonine1}) is reduced to  the definition of the fractional derivative   (\ref{GruLet_der}) without any auxiliary function, then Letnikov concluded that the following alternative holds: either 1) there exists no auxiliary function of the above form, or 2) formula (\ref{Sonine1}) can not be taken  as the general definition of fractional derivative of arbitrary order. Instead, he said that his definition is free of necessity to add an auxiliary function (in spite the fact that this definition is reduced to different form   (\ref{GruLet_int}) and  (\ref{GruLet_der}) whenever $\alpha$ is negative or positive, respectively).

\subsection{Sonine kernel and Sonine integral equations}
\label{Kernel}

In one of his pioneering articles \cite{Abe26} Abel presented solution to the integral equation
\begin{equation}
\label{Abel1}
\int\limits_{0}^{x} \frac{\varphi(\tau) d \tau}{(x - \tau)^{1 - \alpha}} = F(x), \; 0 < \alpha < 1.
\end{equation}
The main component of Abel's method was the following identity
\begin{equation}
\label{Abel2}
\int\limits_{0}^{x} f(t) dt = \frac{\sin \pi \alpha}{\pi} \int\limits_{0}^{x} \frac{d t}{(x - t)^{1-\alpha}}  \int\limits_{0}^{t} \frac{f(\tau) d \tau}{(t - \tau)^{\alpha}},
\end{equation}
where $f(x) = F^{\prime}(x)$.

Sonine tried to generalize Abel's identity (\ref{Abel2}) in order to solve more general equation than integral equation (\ref{Abel1}). He looked for a pair of functions $\sigma(x), \psi(x)$ satisfying the identity
\begin{equation}
\label{SonAbe2}
\int\limits_{0}^{x} f(t) dt = \int\limits_{0}^{x} \psi(x - t) dt \int\limits_{0}^{t} f(\tau) \sigma(t - \tau) d \tau,
\end{equation}
or, what is the same, a pair of functions generating integral representation of unity
\begin{equation}
\label{SonAbe2u}
1 = \int\limits_{0}^{x} \psi(x - t) dt \int\limits_{0}^{t} \sigma(t - \tau) d \tau.
\end{equation}
Sonine described in \cite{Son84} a possible form of the functions  $\sigma(x), \psi(x)$
$$
\sigma(t) = \frac{t^{-p}}{\Gamma(1 - p)} \sum\limits_{k=0}^{\infty} a_k t^k,\;\;\; \psi(t) = \frac{t^{-q}}{\Gamma(1 - q)} \sum\limits_{k=0}^{\infty} b_k t^k,
$$
where $p + q = 1$ and coefficients $a_k, b_k$ are defined by the following relations
$$
a_0 b_0 = 1, \; \sum\limits_{k=0}^{n} \Gamma(k + p) \Gamma(q + n - k) a_{n -k} b_k = 0,\; n = 1, 2, \ldots.
$$
He also applied relation (\ref{SonAbe2u}) for representation of solution to the first kind integral equation with one of these functions as the kernel:
\begin{equation}
\label{SonIntEqu}
\int\limits_{0}^{x} \sigma(x - \tau) \varphi(\tau)  d \tau = f(x).
\end{equation}
Both functions   $\sigma(x), \psi(x)$ are known as Sonine kernels, and integral equation (\ref{SonIntEqu}), generalizing Abel integral equation (\ref{Abel1}), is called Sonine integral equation.  In modern language (see, e.g.  \cite{SamCar03}) a locally integrable function $\sigma(x)$ is called the Sonine kernel if there exists another locally integrable function $\psi(x)$ such that the following identity holds
\begin{equation}
\label{Son_kernel}
\int\limits_{0}^{x} \sigma(x - \tau) \psi(\tau)  d \tau = 1, \; x > 0.
\end{equation}
In fact, the function  $\psi(x)$ is also called the Sonine kernel (sometimes these functions are called the associated Sonine kernels).

Several special examples of Sonine kernel are presented in \cite{Luc21}. We can mention also the article \cite{SamCar03} in which the properties of the Sonine kernel are discussed in modern setting. Several difficulties which one has to overcome by formal application of Sonine's approach to the solution of the corresponding integral equations are discovered and possible ways to overcome these difficulties are shown.

In \cite{Luc21a}, it is introduced the general fractional integrals and derivatives of arbitrary
 order and study some of their basic properties and particular cases. First, a suitable generalization
 of the Sonine condition is presented and some important classes of the kernels that satisfy this
 condition are introduced. In the introduction of the general fractional integrals and derivatives the author follows
 recent approach by Kochubei \cite{Koc11}. The general fractional integrals and derivatives with Sonine kernel  are defined in the Riemann-Liouville form (see \cite{Luc21,Luc21a} and references therein)
 \begin{equation}
 \label{Son_FI}
 \left({\mathbb I}_{\sigma} f\right)(x) = \int\limits_{0}^{x} \sigma(x - t) f(t) d t,\; x > 0,
 \end{equation}
\begin{equation}
 \label{Son_FD}
 \left({\mathbb D}_{\psi} f\right)(x) = \frac{d}{d x} \int\limits_{0}^{x} \psi(x - t) f(t) d t,\; x > 0,
 \end{equation}
 where the functions   $\sigma(x), \psi(x)$ are associated Sonine kernels. Operators (\ref{Son_FI}), (\ref{Son_FD}) are discussed in \cite{Koc11,Luc21,Luc21a} under different conditions on  Sonine kernels and constructions similar not only to Riemann-Liouville type fractional integrals and derivatives, but also to  Dzhrbashian-Caputo type and to Marshaud type.

\subsection{Higher order hypergeometric functions}
\label{HyperF}

The main research interest by Sonine was  the study of the properties of several classes of special functions. His results
served as an impetus for the development of the theory of cylindrical functions (or Bessel-type functions) in the second half of XIX century.
These results are based on the achievements by C. Neumann, O. Schl\"omilch, E. Lommel, H. Hankel, N. Nielsen, L. Schlafli, L. Gegenbauer and others (see e.g. \cite{Kro67}, \cite{NIST}). Sonine defined in \cite{Son80} the cylindrical functions $S_\nu(z)$ as a partial solutions to the following system of functional-differential equations
\begin{equation}
\label{Son_cyl}
\left\{\begin{array}{l}
S_{\nu+1}(z) + 2 S^{\prime}_{\nu}(z) - S_{\nu - 1}(z) = 0, \\
2 \nu  S_{\nu}(z) = z \left[S_{\nu-1}(z) + S_{\nu+1}(z)\right], \\
S_1(z) = - S_0^{\prime}(z),
\end{array}
\right.
\end{equation}
where $z$ is the complex variable and $\nu$ is an arbitrary complex parameter. Sonine proved that these partial solutions admit an integral
representation
\begin{equation}
\label{SonIntRepr}
S_{\nu}(z) = \frac{1}{2\pi i} \int\limits_{a}^{b} \exp\left\{\frac{z}{2}\left(t - \frac{1}{t}\right)\right\} \frac{dt}{t^{\nu + 1}}.
\end{equation}
He found four possible cases for the limits of integration, namely: 1) $\infty \cdot\alpha$, $\infty \cdot \beta$; 2) $- \frac{0}{\alpha}$, $- \frac{0}{\beta}$; 3) $- \frac{0}{\alpha}$, $\infty \cdot \beta$; 4) ${\mathrm{Im}} (z a) = \pm \infty$, ${\mathrm{Im}} (z b) = \pm \infty$, where
in the cases 1)-3) ${\mathrm{Re}} (z a) < 0$,  ${\mathrm{Re}} (z b) < 0$, but in the case 4)  ${\mathrm{Re}} (\nu) > 0$. Sonine denoted the functions obtained in these four cases by $S_{\nu}^{(k)}(z)$ and showed that
$$
S_{\nu}^{(1)}(z) = J_{\nu}(z),\; S_{\nu}^{(2)}(z) = e^{-\nu \pi i} J_{-\nu}(z),\; S_{\nu}^{(3)}(z) = \frac{1}{2} H_{\nu}^{(1)}(z),\; S_{\nu}^{(4)}(z) = J_{\nu}(z).
$$
The above integral representation (\ref{SonIntRepr}) is called Sonine integral representation which is a source for obtaining new representations for cylindrical functions (see \cite{NIST}) as well as for calculation of certain definite integrals. Among these integrals are those known as the first and the second Sonine integrals, respectively (or classical Sonine formulas), see e.g. \cite{Gra19}:
\begin{equation}
\label{SonIntegral1}
J_{\nu + \mu + 1}(a q) = \frac{q^{\nu}}{2^{\nu} \Gamma(\nu + 1) a^{\nu + \mu + 1}} \int\limits_{0}^{a} J_{\mu}(q x) (a^2 - x^2)^{\nu} x^{\mu + 1} dx,
\end{equation}
\begin{equation}
\label{SonIntegral2}
 \int\limits_{0}^{a} J_{\mu}(q x) J_{\nu}[z\sqrt{a^2 - x^2}] (a^2 - x^2)^{\frac{\nu}{2}} x^{\mu + 1} dx = a^{\nu + \mu + 1} q^{\mu} z^{\nu}
 \frac{J_{\nu + \mu + 1}(a \sqrt{q^2 + z^2})}{(\sqrt{q^2 + z^2})^{\nu + \mu + 1}},
\end{equation}
where ${\mathrm{Re}} \nu, {\mathrm{Re}} \mu > - 1$. Sonine formulas find interest in different questions of analysis (e.g. in Dunkl theory, see \cite{RoeVoi20}, or in the study of Levy processes \cite{BN-M-R}).

There exist several multivariate extensions of the classical Sonine
integral representation for Bessel functions of some index $\mu + \nu$ with respect to
such functions of lower index $\mu$ (see, e.g. \cite{RoeVoi18}). For Bessel functions on matrix cones, Sonine
formulas involve beta densities $\beta_{\mu,nu}$ on the cone and trace already back to
Herz.

Several important results dealing with properties of $\Gamma$-function were obtained by Sonine during his career. They are based on the study of the solution to the difference equation
\begin{equation}
\label{Son_diff}
F(x + 1) - F(x) = f(x).
\end{equation}
In these works Sonine followed the idea by Binet (1838) who examined the relation
$$
\log \Gamma(x + 1) - \log \Gamma(x) = \log x.
$$

Sonine found \cite{Son81}, in particular, the form of the remainder factor in the product representation of $\Gamma$-function
\begin{equation}
\label{Son_prod}
\Gamma(x + 1) = \frac{n! (n+1)^x}{(x+1)(x+2)\cdots (x+n)} \frac{\left(1 + \frac{x}{n+1}\right)^{x+n+\theta}}{\left(1 + \frac{1}{n+1}\right)^{x(1+n+\theta)}}, \; x\in {\mathbb R}, \; 0 < \theta < 1.
\end{equation}
In his article on Bernoulli polynomials Sonine obtained one more representation,  related to $\Gamma$-functions (this formula was rediscovered by Ch.Hermite in 1895)
\begin{equation}
\label{Son_sum}
\log \frac{\Gamma(x + y)}{\Gamma(y)} = x \log y + \sum\limits_{k=2}^{n} \frac{(-1)^k \varphi_k(x)}{(k - 1) k y^{k  -1}} + R_n(x, y),
\end{equation}
where $\varphi_k(x)$ are Bernoulli polynomials defined by Sonine using difference equation
$$
\varphi_k(x + 1) - \varphi_k(x) = k x^{k - 1},\; \varphi_k(0) = 0, k = 1, 2, \ldots
$$

The book \cite{Son54} contains a number of the most important articles by N. Ya Sonine as well as a survey on his other research.

\section{Nekrasov's contribution to Fractional Calculus}
\label{Nekrasov}

In \cite{Nek88} Nekrasov proposed a new definition of the general differentiation. In fact, this definition include as a special case
Letnikov's definition. Main idea by Nekrasov is to define the derivative by using integration along closed contour $L$ crossing the point $x$ and encircled a group of singular points of the differentiable function $f(x)$. This definition gives, in fact, the differentiation with respect to a doubly connected domain, which is free of the singular points of $f(x)$, and contains above said contour $L$. Therefore, Nekrasov used the ideas by Sonine (to take into account the properties of the analytic continuation of the given function and to apply the properties of functions in complex domains). The main aim of Nekrasov's construction is to extend the class of functions to which the general differentiation can be applied.

It should be noted that the construction proposed by Nekrasov is fairly complicated and need to use properties of the functions on Riemann surfaces. It follows from the properties of functions to which Nekrasov tried to apply his definition. Starting point of his construction is the notion of classes $(q, \mu)$ of function. Let $L$ be a closed contour encircled a group of singular points of the function $f(z)$. Let the function $f(z)$ has the following property: if the point $z$ makes a complete  detour along $L$ in counter clockwise direction, then the function $f(z)$, continuously changing, gains the multiplier $e^{2\pi q i}$. Then this function is of class $(q, 0)$. Thus, any function of the class  $(q, 0)$ can be represented in the form $f(z) = (z - a)^{q} \phi(z)$, where $a$ lies inside $L$ and $\phi(z)$ is of the class  $(0, 0)$. The function of the form $f(z) = (z - a)^{q} \log^{\mu} z \phi(z)$, with $\phi(z)$ being of the class  $(0, 0)$, is said to belong to the class $(q, \mu)$ (with $q$ being the power index, and $\mu$ being the logarithmic index which is supposed to be nonnegative integer number). It is clear that if the function $f(z)$ belongs to the class  $(q, \mu)$, then it belongs to any class  $(q \pm m, \mu)$, $m\in {\mathbb N}$. Clearly, this definition depends on the choice of the contour $L$.

The function $f(z)$, which can be represented in form of a sum of finite number $n$ of functions belonging to different classes with respect to the contour $L$, is said to be reducible to $n$ classes (or simply reducible).

Let the function $f(z)$ be reducible to $n$ classes with zero logarithmic indices, i.e.
\begin{equation}
\label{Nek01}
f(z) = f_0(z) + f_1(z) + \ldots + f_{n-1}(z),
\end{equation}
where $f_j(z)$ is of class $(q_j, 0)$. Then we have the following representation
\begin{equation}
\label{Nek1}
f(z) + \int\limits_{(L^k)}^{(z)} \frac{d f(t)}{d t} dt = \alpha_{0}^{k} f_0(z) + \alpha_{1}^{k} f_1(z) + \ldots + \alpha_{n-1}^{k} f_{n-1}(z),
\end{equation}
where integration is performed along the contour $L$, traversable  $k$-times in counter clockwise direction starting from the point $z$, $\alpha_{j} = e^{2\pi q_j i}$. By assigning the values $0, 1, \ldots, n - 1$ to the parameter $k$ we obtain the system of equations sufficient for determination of functions $f_0(z), f_1(z), \ldots, f_{n-1}(z)$.

Let the function  $f(z)$ be reducible to $n$ classes with non-zero logarithmic indices, i.e.
\begin{equation}
\label{Nek02}
f(z) = \sum\limits_{s=0}^{n_0 - 1} (z - a)^{q_{s,0}} \phi_{s,0}(z) +  \sum\limits_{s=0}^{n_1 - 1} (z - a)^{q_{s,1}} \log (z - a) \phi_{s,1}(z) + \ldots +  \sum\limits_{s=0}^{n_\mu - 1} (z - a)^{q_{s,\mu}} \log^{\mu} (z - a) \phi_{s,\mu}(z),
\end{equation}
where $n_0 + n_1 + \ldots + n_\mu = n$ and all functions $\phi_{s,j}(z)$ are of class $(0, 0)$. Then we have the following representation
\begin{equation}
\label{Nek2}
f(z) + \int\limits_{(L^k)}^{(z)} \frac{d f(t)}{d t} dt =  \sum\limits_{s=0}^{n_0 - 1} \alpha_{s,0}^{k} (z - a)^{q_{s,0}} \phi_{s,0}(z) + \sum\limits_{s=0}^{n_1 - 1} \alpha_{s,1}^{k} (z - a)^{q_{s,1}} \{2k \pi i + \log (z - a)\} \phi_{s,1}(z) + \ldots
\end{equation}
$$
+  \sum\limits_{s=0}^{n_\mu - 1}  \alpha_{s,\mu}^{k} (z - a)^{q_{s,\mu}} \{2k \pi i + \log (z - a)\}^{\mu} \phi_{s,\mu}(z),
$$
where integration is performed along the contour $L$, traversable  $k$-times in counter clockwise direction starting from the point $z$, $\alpha_{s,j} = e^{2\pi q_{s,j} i}$. By assigning the values $0, 1, \ldots, n - 1$ to the parameter $k$ we obtain the system of equations sufficient for determination of functions $\phi_{s,0}(z), \phi_{s,1}(z), \ldots, \phi_{s,n-1}(z)$.

Therefore, in both cases we have a finite sum representation of the function $f(z)$ of the considered form. Now the question is to define the integral/derivative of arbitrary order of each components of the representation (\ref{Nek01}) or (\ref{Nek02}). Moreover, any function $f(z)$ of the class $(q, \mu)$ can be determined as the following limit
\begin{equation}
\label{Nek_lim}
f(z) = \lim\limits_{h\rightarrow 0} \left[(z - a)^q \left(\frac{(z -a)^h - 1}{h}\right)^{\mu}\right].
\end{equation}
Thus we have the limit of the finite sum of functions belonging to the classes $(q, 0), (q + h, 0), \ldots, (q + \mu h, 0)$. Therefore the definition of the integral/derivative of arbitrary order of the reducible functions can be completely described if one can define the definition of a function of class $(q, 0)$. {Nekrasov noted that his construction of his derivative generally speaking can not be rigorously defined in the case when $f(z)$ is reducible to infinite number of classes.}

For this definition Nekrasov used Letnikov's formulas. The only difference is that the contour of integration is now specially deformed curve on the Riemann surfaces (which depend on the order of the derivative, i.e. can be either finite-sheeted or infinite-sheeted).

In the article \cite{Nek89} Nekrasov apply his construction to determination of the solution of the following differential equation
\begin{equation}
\label{Nek_DE}
\sum\limits_{s=0}^{n} (a_s + b_s x) x^s \frac{d^s y}{d x^s} = 0,
\end{equation}
which is highly related to the generalized hypegeometric function ${}_p F_{q}(z)$.

\vspace{5mm}
\noindent
{\Large\bf{Short biographies}}

\subsection*{Letnikov Alexey Vasil'evich (1837-1888), Russian mathematician. A short biography.}

A.V. Letnikov was born on January 1st, 1837, in Moscow, Russia. When Alexey was 8 years old his father died. His mother tried to give education to Alexey and his sister. Mother sent Alexey to grammar school in 1847. In spite of his evident abilities he was not too successful in education. Therefore he was moved to Konstantin's land-surveyors institute (full-time provisional military type institute). That was a second rank educational establishment. The director of the discovered high interests of Alexey to mathematics and supported his growth in the subject. The director decided to prepare him to career of teacher in mathematics in Konstantin's land-surveyors institute. To get the corresponding position Letnikov was sent to Moscow University and study mathematics during two years (1856-1858) there as an extern student.

After graduation he was sent to Paris in order to extend his knowledge at the most known mathematical center during two years and to study the structure and the content of the technical education in France. In Paris, Letnikov took lectures of many well-known mathematicians (Liouville among them) in the Ecole Polytechnique, College de France and Sorbonne.

Returning from Paris in December 1860 he was appointed as a teacher in the engineering class of the Konstantin's land-surveyors institute and started to teach Probability Theory. Letnikov actively participated in mathematical life in Moscow. In particular, he was among the founders of Moscow Mathematical Society in 1864. In 1863 it was approved a new Statute of Higher Education. Among other regulations it supposed to enlarge a number of chairs at universities and to recruit new university teachers. To get a position at university one needs either to pass graduation gymnasium's exams or to receive the degree at a foreign university. Letnikov decided to use the second possibility. In 1867 he defended PhD "\"{U}ber die Bedingungen der Integrabilit\"{a}t einiger Differential-Gleichungen" at Leipzig University. In 1868 he got a position at recently reopened Imperial Technical College (now Bauman's Technical University). Letnikov was working at this College up to 1883 when he moved to Alexandrov's Commercial College sharing this job with a part-time teacher at Konstantin's land-surveyors institute and at the Imperial Technical College.

It was active time for him, he was awarded the degree of a state councillor, got the order of Saint Stanislav and was appointed in 1884 as a corresponding member of St.-Petersburg Academy of Sciences (by recommendation of V.Imshenetsky, V.Bunyakovsky and O.Backlund). At the end of 1880s Letnikov should receive a state pension and supposed to leave teaching and to concentrate on research. He was dreaming to get a position at Moscow University. It was not happen since at the opening ceremony of a new building of Alexandrov's Commercial College he caught a cold. He had no serious illness before and continued to deliver lectures. But this time the illness was strong enough and he died on February 27th, 1888.

\subsection*{Sonine Nikolai Yakovlevich (1849-1915), Russian mathematician. A short biography.}

N.Ya.Sonine was born on February 10th, 1849, in Tula, Russia.

Studied at Physical-mathematical Faculty of Moscow University (1865-1869). After graduation continued research in Moscow University during two years and in 1871 defended Master Thesis "On expansion of function into infinite series". In June 1871 he became Associate Professor (dozent) of Warsaw University.

In 1873 he was sent to Paris to continue research study. In Paris he got lectures by Liouville, Hermite, Bertrand, Serre and Darboux. In September 1874 in Moscow University N.Ya. Sonine defended PhD Thesis "On integration of partial differential equations of the second order". In 1877 he became extra-ordinary Professor of Warsaw University and since 1879 he was an ordinary Professor of Warsaw University. In 1891 he resigned from his position at Warsaw University, but still continued his scientific work. In 1891 N.Ya. Sonine was elected a corresponding member of Academy of Sciences and 1893 he became an academician of St.-Petersburg Academy of Sciences (by recommendation of P.L.Chebyshev). In 1890 he was awarded by V.Ya.Buniakovsky Prize for the Best Results in Mathematics.

Since 1899 N.Ya. Sonine occupied different administrative positions mostly in education. He died in February 18th, 1915, in St.-Petersburg.

\subsection*{Nekrasov Pavel Alekseevich (1853-1924), Russian mathematician and philosopher. A short biography.}

P.A.Nekrasov was born on February  1st (13th), 1853, in Ryazan region, Russia.

After graduation at Ryazan Orthodox seminary 1874 he entered Physical-mathematical faculty of Moscow University. In 1878 he graduated at Moscow University with degree of candidate of sciences and was left at the department of pure mathematics for preparation to professorship. From August 1879 P.A.Nekrasov shared his research with teaching mathematics at the private Voskresensky's real school. In 1883 he defended master thesis "Study of the equation $u^m - p u^n - q = 0$". For this work he was awarded by V.Ya.Buniakovsky Prize for the Best Results in Mathematics.

In 1985 P.A.Nekrasov became a Privatdozent at Moscow University (having defended his Russian PhD "On Lagrange series" in 1886) and, in 1886, he got a position of an associate professor (extraordinary professor) at Moscow University. In 1890 he received a full professorship. In 1893 he became the rector of Moscow University. After his term as the rector, he actually wanted to retire, but was not allowed to. He also taught 1885-1891 Probability Theory and Higher Mathematics at the Land-surveyors institute.

From 1898 he was almost only with administrative duties for the Ministry of Education (he was curator of the Moscow University and responsible for the schools in Moscow and the surrounding area) and moved in 1905 to Saint Petersburg as a member of the Council of the Ministry of Education. After the Russian Revolution, he tries to adapt to the new realities, dealt with mathematical economics (which he lectured in 1918-19) and studied Marxism. He died of pneumonia 20th December, 1924, in Moscow.


{\bf Acknowledgments}

The work has been supported by Belarusian Fund for Fundamental
 Scientific Research  through the grant F20R-083.





{\bf Conflicts of Interests}

The author declares no conflict of interest.

\bibliographystyle{mdpi}
\makeatletter
\renewcommand\@biblabel[1]{#1. }
\makeatother


\begin{thebibliography}{999} 

\bibitem[SKM]{SKM} Samko, S. G., Kilbas, A. A. and Marichev, O. I. Fractional Integrals and Derivatives. Theory and Applications, Gordon and Breach Science Publishers, Yverdon, 1993, translated from Russian edition, Fractional Integrals and Derivatives and Some of Their Applications, Nauka i Tekhnika, Minsk, 1987.


\bibitem[Let68a]{Let68a} Letnikov, A.V. On historical development of the theory of differentiation of an arbitrary order, Mat. Sb., 3 (2), 85--112 (1868) (in Russian)


\bibitem[Pot03]{Pot03} Potapov, A. A. A short essay on the origin and formation of the theory of fractional integro-differentiation, Nonlinear World, 1 (1-2), 69--81 (2003) (in Russian)

\bibitem[Pot12]{Pot12} Potapov, A. A. Essays on the development of fractional calculus in the A.V. Letnikov's works. To 175th anniversary of A.V.Letnikov, Radioelectronic, Nanosystems, Information Technologies (Radioelektronika, Nanosistemy, Informacionnye Tehnologii), 4 (1), 3--102 (2012) (in Russian)


\bibitem[LetChe11]{LetChe11}
 Letnikov, A.V.; Chernykh, V.A. {The Foundation of Fractional Calculus (with applications to the theory
 of oil and gas production, underground hydrodynamics and dynamics of biological systems)}, Moscow:
 Neftegaz, 2011. - 429 p. (in Russian)

\bibitem[Let68]{Let68} Letnikov, A.V. Theory of differentiation of an arbitrary order, Mat. Sb., 3 (1), 1--68 (1868) (in Russian)

\bibitem[Let73]{Let73} Letnikov A,V. To explanation of the main statements of the theory of differentiation of an arbitrary order, Mat. Sb., 6 (4), 413--445 (1873) (in Russian)


\bibitem[Nek88]{Nek88} Nekrasov, P. A. General differentiation, Mat. Sb., 14 (1), 45--166 (1888) (in Russian)


\bibitem[Son72]{Son72} Sonine, N. Ya. On differentiation of arbitrary order, Mat. Sb., 6 (1), 1--38 (1872) (in Russian)


\bibitem[Vas-Zak61]{Vas-Zak61} Vashchenko-Zakharchenko, M.E. On fractional differentiation, Quarterly J.  Pure and Appl. Math., Ser. 1, IV, 237--243 (1861)

\bibitem[Gru67]{Gru67} Gr\"unwald, A.K. Uber "begrentze" Derivationen und deren Anwendung, Zeitschrift fur Angew. Math. Phys., 12, 441--480 (1867)


\bibitem[GKMR2]{GKMR2} Gorenflo R., Kilbas A.A., Mainardi F., Rogosin S. Mittag-Leffler Functions: Related Topics and Applications.   2nd extended and revised edition.  Berlin - New York: Springer, 2020. - 540 p.


\bibitem[Lio32a]{Lio32a} Liouville, J. Memoire sur quelques Questions de
Geometrie et de Mecanique, et sur un nouveau genre de Calcul
pour resoudre ces Questions. J. Ecole Polytech., 13 (21), 1--69 (1832)

\bibitem[Lio32b]{Lio32b} Liouville, J. Memoire sur le Calcul des differentielles a
indices quelconques. J. Ecole Polytech, 13 (21), 71--162 (1832)

\bibitem[Lio34]{Lio34} Liouville, J. Memoire sur une formule d analyse. J. Reine
Angew. Math. (Grelle's Journal),  12, 273--287 (1834)


\bibitem[NekPok89]{NekPok89} Nekrasov, P. A., Pokrovskii, P. M. On examination of manuscripts by A.V. Letnikov, presented after his death to the Moscow Mathematical Society, Mat. Sb., 14 (2), 202--204 (1889) (in Russian)


\bibitem[Let89a]{Let89a} Letnikov, A.V. On integration of the equation $(a_n + b_n x) \frac{d^n y}{d x^n} + (a_{n-1} + b_{n-1} x) \frac{d^{n-1} y}{d x^{n-1}} + \ldots + (a_0 + b_0 x) y = 0$, Mat. Sb., 14 (2), 205--215 (1889) (in Russian)

\bibitem[Abe26]{Abe26}  Abel, N. H. Aufl\"osung einer mechanischen Aufgabe, Journal f\"ur die reine und angewandte Mathematik, 1, 153--157 (1826)


\bibitem[Son84]{Son84} Sonine, N. Ya. Sur la g\'{e}n\'{e}ralisation d'une formule d'Abel, Acta Math., 4, 171--176 (1884) (in French)


\bibitem[SamCar03]{SamCar03} Samko S.G, Cardoso, R.P. Integral equations of the first kind of Sonine type, IJMMS, 57, 3609--3632 (2003)


\bibitem[Luc21]{Luc21} Luchko Y. General Fractional Integrals and Derivatives with the Sonine Kernels. Mathematics, 9 (6), 594 (2021)

\bibitem[Luc21a]{Luc21a} Luchko, Yu. General Fractional Integrals and Derivatives of Arbitrary Order. Symmetry, 1,
0  (2021).

\bibitem[Koc11]{Koc11}  Kochubei, A. N. General fractional calculus, evolution equations, and renewal processes. Integr. Equa. Operator Theory, 71,
583--600  (2011)


\bibitem[Kro67]{Kro67} Kropotov, A. I. Nikolai Yakovlevich Sonine. Leningrad: Nauka, 1967. - 136 p. (in Russian)


\bibitem[NIST]{NIST} 
{\it NIST Handbook of Mathematical Functions}. Edited by Frank
W.J. Olver (editor-in-chief), D.W. Lozier, R.F. Boisvert, and C.W.
Clark. Gaithersburg, Maryland, National Institute of Standards and
Technology, and New York, Cambridge University Press, 951 + xv
pages and a CD, (2010)


\bibitem[Son80]{Son80} Sonine, N. Ya. Recherches sur les fonctions celindriques et le d\'{e}veloppment des fonctions continues en series, Math. Ann., Bd. 16, 1--80 (1880)


\bibitem[Gra19]{Gra19} Grandits, P. Some notes on Sonine-Gegenbauer integrals, Int. Transf. Spec. Funct., 30 (2), 128--137 (2019)


\bibitem[RoeVoi20]{RoeVoi20} R\"osler M., Voit, M. Sonine Formulas and Intertwining Operators in Dunkl Theory, International Mathematics Research Notices, rnz313 (2020), arXiv:1902.02821v2 [math.CA] 22 Oct 2019


\bibitem[BN-M-R]{BN-M-R} Barndorff-Nielsen, O. E., Mikosch, T., Resnick, S. I. (Eds.) Levy processes. Theory and applications. Boston
(MA): Birkh\"auser, 2001. - 418 p.


\bibitem[RoeVoi18]{RoeVoi18} R\"osler M., Voit, M. Beta distributions and sonine integrals for Bessel functions on symmetric cones, Studies in Appl. Math., 141 (4), 474--500 (2018)


\bibitem[Son81]{Son81} Sonine, N. Ya. Note sur une formule de Gauss, Bull. Soc. Math. France, 9, 162--166 (1881)  (in French)



\bibitem[Son54]{Son54} Sonine, N. Ya. Research on cylindric functions and orthogonal polynomials. Edited and commented by N. I. Akhieser. Moscow: GITTL,  1954. - 244 p. (in Russian)

\bibitem[Nek89]{Nek89} Nekrasov, P. A. Application of the general differentiation to the integration of the equation of the type $\sum (a_s + b_s x) x^s D^s y = 0$, Mat. Sb., 14 (3), 344--393 (1889) (in Russian)







\end{thebibliography}


%


%

\end{document}